\documentclass{amsart}
\usepackage{amsfonts, amssymb}
\usepackage{eucal,latexsym}

\newtheorem{thm}{Theorem}[section]
\newtheorem{cor}[thm]{Corollary}
\newtheorem{lem}[thm]{Lemma}

\theoremstyle{definition}
\newtheorem{defn}[thm]{Definition}
\theoremstyle{remark}

\numberwithin{equation}{section}
\newcommand{\st}{{\rm st}}
\newcommand{\CC}{\mathcal C}
\newcommand{\hh}{{\bf H}}
\newcommand{\HH}{\mathcal H}

\newcommand{\RR}{\mathbb R}

\begin{document}

\title[Compression of uniform embeddings into Hilbert space]
{Compression of uniform embeddings into Hilbert space}

\author{N.~Brodskiy}
\address{Department of Mathematics, University of Tennessee, Knoxville, TN 37996, USA}
\email{brodskiy@math.utk.edu}

\author{D.~Sonkin}
\address{Department of Mathematics, University of Virginia, Charlottesville, VA 22903, USA}
\email{ds5nd@virginia.edu}
\thanks{The work of the second named author was supported in part by the NSF grant 0245600 of A.Yu. Ol'shanskii and M.V. Sapir.}




\begin{abstract}
If one tries to embed a metric space uniformly in Hilbert space,
how close to quasi-isometric could the embedding be? We answer
this question for finite dimensional CAT(0) cube complexes and for
hyperbolic groups. In particular, we show that the Hilbert space
compression of any hyperbolic group is 1.
\end{abstract}

\maketitle


\section{Introduction}


The study of uniform embeddings of metric spaces into a Hilbert
space was originated by Gromov~\cite{gro}.

\begin{defn}
Let $(X,d)$ be a metric space and $\HH$ be a Hilbert space. A map
$f\colon X \to \HH$ is said to be a {\it uniform embedding} if
there exist non-decreasing functions $\rho,\delta\colon\RR_+ \to
\RR_+$ such that
\begin{itemize}
\item [(1)] $\rho(d_X(x,y)) \le d_\HH(f(x),f(y)) \le
\delta(d_X(x,y))$ for all $x,y \in X$;

\item [(2)] $\lim\limits_{t \to +\infty} \rho(t)=+\infty$.

\end{itemize}
\end{defn}

Gromov asked whether a finitely generated group (viewed as a
metric space with a word length metric) that is uniformly
embeddable into a Hilbert space satisfies the Novikov Conjecture.
Yu ~\cite{yu} answered this question affirmatively, showing that
for a finitely generated group equipped with a word length metric,
uniform embeddability into a Hilbert space implies the Coarse
Baum-Connes Conjecture and the Novikov Conjecture. Uniform
embeddability of finitely generated groups into a Hilbert space
has been studied extensively since then (see~\cite{DadGue,
gue-kam, ags, DadGue1, chen-dad-gue-yu, cam-nib} and references
therein).

The definition above suggests that there exist two useful
real-valued functions associated to the map $f$: its {\it
dilatation}
$$\delta_f(t)=\sup\{d_\HH(f(x),f(y))\mid d_X(x,y)\le t\}$$
and its {\it compression}
$$\rho_f(t)=\inf\{d_\HH(f(x),f(y))\mid d_X(x,y)\ge t\}.$$

If $X$ is a quasi-geodesic metric space (in particular, a finitely
generated group with a word metric), then the dilatation of its
uniform embedding in a Hilbert space is dominated by a linear
function~\cite{gro, gue-kam}. We consider the following question:
how close to a linear function can the compression of a uniform
embedding of a group into a Hilbert space be? In other words, how
close to a quasi-isometric embedding can a uniform embedding be?

We are interested in asymptotic behavior of the compression
function. Therefore it is convenient to introduce a (partial)
relation on functions $g,h\colon\RR_+ \to \RR_+$ as follows: we
write $f\preceq g$ if there exist numbers $C$ and $M$ such that
$f(t)\le Cg(t)$ for all $t\ge M$. We write $f\sim g$ if $f\preceq
g\preceq f$.

Guentner and Kaminker~\cite{gue-kam} introduced a geometric
invariant $R(G)$ of a finitely generated group $G$ called Hilbert
space compression. $R(G)$ is the supremum of all numbers
$\alpha\ge 0$ for which there exists a uniform embedding $f\colon
G \to \HH$ with the compression $\rho_f(t)\succeq t^\alpha$. This
invariant is a number between 0 and 1 and it parameterizes the
difference between $G$ being uniformly embeddable in a Hilbert
space and being exact \cite{gue-kam}.

It was shown in~\cite{gue-kam} is $R(F_2)=1$ where $F_2$ is a free
group on two generators. The authors proved it by providing a
sequence of uniform embeddings $f_n\colon F_2 \to \HH$ with
$\rho_{f_n}(t)\succeq t^{1-\frac{1}{n}}$. A question remained of
whether there is one uniform embedding $f\colon F_2 \to \HH$ with
$\rho_{f}(t)\succeq t^{1-\frac{1}{n}}$ for all $n$. We construct
such an embedding in Section~\ref{trees} by modifying the
construction from~\cite{gue-kam}. Moreover we estimate the
compression function of our embedding as
$$\rho_f(t)\succeq\frac{t}{\sqrt{\ln t}\cdot\ln\ln t}$$
It follows from results of Bourgain~\cite{bourg} that for any
uniform embedding $f\colon F_2 \to \HH$ we have
$$\rho_f(t)\preceq\frac{t}{\sqrt{\ln t}}.$$
It is not known whether there exists a uniform embedding of $F_2$
in $\HH$ with the compression function $\sim \frac{t}{\sqrt{\ln
t}}$.

We show in Section~\ref{HypGroups} that any hyperbolic group
embeds uniformly in a Hilbert space with the compression function
$\rho_f(t)\succeq\frac{t}{\sqrt{\ln t}\cdot\ln\ln t}$. In
particular, the Hilbert space compression of any hyperbolic group
is equal to 1. Notice that for any uniform embedding of a
non-elementary hyperbolic group in a Hilbert space we have
$\rho_f(t)\preceq\frac{t}{\sqrt{\ln t}}$.

Finally, Section~\ref{CAT(0)} is devoted to uniform embeddings of
groups acting on CAT(0) cubical complexes. It was proven by
Campbell and Niblo~\cite{cam-nib} that any discrete group $G$
acting properly, co-compactly on a finite dimensional CAT(0)
cubical complex has Hilbert space compression $R(G)=1$. They
provided a sequence $f_n$ of uniform embeddings with compressions
$\rho_{f_n}(t)\succeq t^{1-\frac{1}{n}}$. We modify their
construction to embed such a group uniformly in a Hilbert space
with compression $\rho_f(t)\succeq\frac{t}{\sqrt{\ln t}\cdot\ln\ln
t}$. It is proved recently by Sageev and Wise~\cite{sag-wise} that
groups acting properly and co-compactly on finite dimensional
CAT(0) cubical complexes satisfy the Tits alternative provided the
orders of finite subgroups are uniformly bounded. We note that for
such a group $G$ either a quasi-isometric embedding in a Hilbert
space exists or any uniform embedding of $G$ in Hilbert space has
compression $\rho_f(t)\preceq\frac{t}{\sqrt{\ln t}}$.

The authors wish to thank Mark Sapir, Graham Niblo, Alexander
Ol'shanskii, Denis Osin, Michah Sageev and Sergey Borodachov for
interesting conversations during the course of this research.

\section{Uniform embeddings of trees}\label{trees}

All trees considered in this section are locally finite, with
geodesic metric and with all edges of length 1.

We start with an observation (for the proof see, for example,
\cite[Proposition~2.9]{gue-kam}).

\begin{lem} \label{emb}
Let $X$ and $Y$ be metric spaces, and assume that $X$ is geodesic.
If $f\colon X\to Y$ is a uniform embedding, then
$\delta_f(t)\preceq C\cdot t$ for some constant $C>0$.
\end{lem}

Results of Bourgain~\cite{bourg} (see~\cite{LinSaks} for the
theorem we use and its short proof) imply the following

\begin{thm} \label{Bourgain}
For any uniform embedding $f$ of a complete infinite binary tree
into a Hilbert space,
$$
\rho_f(t)\preceq\frac{\delta_f(t)}{\sqrt{\ln t}}.
$$
\end{thm}

Since complete infinite binary tree embeds isometrically into the
Cayley graph of the free group $F_2$ relative to the standard
generating set, in view of Lemma~\ref{emb} we have

\begin{cor} \label{compression tree}
If $f\colon F_2\to \HH$ is a uniform embedding, then
$\rho_f(t)\preceq\frac{t}{\sqrt{\ln t}}$.
\end{cor}

In particular,

\begin{cor} 
There is no quasi-isometric embedding of the free group $F_2$ into
Hilbert space.
\end{cor}

We need the following technical lemma.

\begin{lem} \label{calc}
Let $\xi(t)=\frac{\sqrt{t}}{\sqrt{\ln t}\cdot\ln\ln t}$. Then
\begin{itemize}
\item [(1)] $\sum_{j=2}^{\infty}(\xi(j+1)-\xi(j))^2<\infty$;

\item [(2)] there exists a constant $C$ such that for any large
enough integer $N$,
$$\sum_{i=1}^{N}\xi^2(i)\ge\frac{1}{2}N\cdot\xi^2(N)-C$$

\end{itemize}
\end{lem}

\begin{proof} (1) It is enough to estimate the sum starting with $j=4$:
$$\sum_{j=4}^{\infty} (\xi(j+1)-\xi(j))^2=\sum_{j=4}^{\infty}
(\int_j^{j+1}\xi'(t)\; dt)^2\le\sum_{j=4}^{\infty}
\int_j^{j+1}[\xi'(t)]^2\; dt=\int_4^\infty[\xi'(t)]^2\; dt$$

In view of the inequality $\xi'(t)\le \frac{\xi(t)}{2t}$ we get

$$\int_4^\infty[\xi'(t)]^2\; dt\le \int_4^\infty\frac{1}{t\cdot \ln t\cdot (\ln\ln
t)^2}\; dt=\int_{\ln\ln 4}^\infty\frac{du}{u^2}<\infty$$

(2) One can check that there exists $M$ such that the function
$\xi(t)$ is increasing for $t\ge M$. Therefore, for any $N \ge M$,
$$\sum_{i=1}^{N} \xi^2(i)\ge\sum_{i=1}^{M}
\xi^2(i)+\int_{M}^{N}\xi^2(t)\; dt$$ The inequality
$\xi^2(t)\ge\frac{1}{2}(t\cdot \xi^2(t))'$ implies that
$$\int_{M}^{N}\xi^2(t)\; dt
\ge\int_{M}^{N}\frac{1}{2}(t\cdot\xi^2(t))'\; dt=\frac{1}{2}N\cdot
\xi^2(N)-\frac{1}{2}M\cdot \xi^2(M)$$
\end{proof}

\begin{thm}\label{embed tree}
Let $T$ be a locally finite tree with geodesic metric and all
edges of length 1. There exists a uniform embedding $f$ of $T$
into Hilbert space $\HH$ with $$\rho_{f}(t) \succeq
\frac{t}{\sqrt{\ln t}\cdot\ln\ln t}$$
\end{thm}

\begin{proof} We are going to define our embedding on vertices of
$T$ and extend it linearly to edges. Let us fix a basis in the
space $\HH$ and enumerate the elements of this basis by all edges
of the tree $T$. If $v$ is an edge of the tree, we denote by
$\vec{e}_v$ the corresponding unit vector of the basis.

Fix a base vertex $O$ of $T$. For any vertex $V$ consider the
geodesic path $[VO]$ joining $V$ and $O$. Let us enumerate the
edges of $[VO]$ by $v_1$, $v_2$, ..., $v_{\|V\|}$ in the order
from $V$ to $O$, where $\|V\|$ the length of $[VO]$. Fix a {\it
weight function} $$ \xi(t)=\frac{\sqrt{t}}{\sqrt{\ln t}\cdot\ln\ln
t} $$ and define the embedding $f$ on the vertex $V$ as
$$ f(V)=\sum_{i=1}^{\|V\|} \xi(i)\cdot \vec{e}_{v_i}. $$

In order to show that $f$ is a uniform embedding we need to
estimate its dilatation $\delta_f$ from above and its compression
$\rho_f$ from below.

First we show that $\delta_f(t)\preceq t$. To prove this we show
that there is a constant $C$ such that the expansion of any edge
of $T$ under the map $f$ does not exceed $C$. Let $U$ and $V$ be
adjacent vertices of the tree $T$. We are going to find a constant
$C$ such that $d(f(U),f(V))\le C\cdot d(U,V)=C$. Recall that we
denote by $v_1$, $v_2$, ..., $v_{\|V\|}$ (respectively by $u_1$,
$u_2$, ..., $u_{\|U\|}$) the edges of the path $[VO]$
(respectively $[UO]$) in the order from $V$ (resp. $U$) to $O$. We
may assume that $\|U\|=\|V\|-1$ and therefore $u_i=v_{i+1}$ for
every $i=1,...,\|U\|$. Since
$$ f(V)-f(U)=\sum_{i=1}^{\|V\|} \xi(i)\cdot \vec{e}_{v_i}-
\sum_{j=1}^{\|U\|} \xi(j)\cdot \vec{e}_{u_j}=\xi(1)\cdot
\vec{e}_{v_1}+\sum_{j=1}^{\|U\|} (\xi(j+1)-\xi(j)) \vec{e}_{u_j},
$$ we have
$$d(f(U),f(V))=\sqrt{\xi^2(1)+\sum_{j=1}^{\|U\|} (\xi(j+1)-\xi(j))^2}$$
and the existence of the constant $C$ follows from part (1) of
Lemma~\ref{calc}.

Now we estimate the function $\rho_f$. Let $U$ and $V$ be any
vertices of the tree $T$ and $S$ be the vertex of $T$ such that
$[UO]\cap [VO]=[SO]$. Then the path from $U$ to $V$ is a union of
paths from $U$ to $S$ and from $S$ to $V$. Without loss of
generality we assume that $|SV|\ge \frac{1}{2}|UV|$.

We estimate $d(f(U),f(V))$ by looking at coordinates corresponding
to edges of the path from $V$ to $S$ (these are edges $v_1$,
$v_2$, ..., $v_{|SV|}$). These coordinates of the point $f(U)$ are
all 0. So,
$$ [d(f(U),f(V))]^2\ge \sum_{i=1}^{|SV|} \xi^2(i). $$
Part (2) of Lemma~\ref{calc} implies
$$\sum_{i=1}^{|SV|} \xi^2(i)\ge \frac{1}{2}|SV|\cdot
\xi^2(|SV|)+{\text{const}}$$ Since the function $\xi^2(t)$
increases for large $t$, we can put $\frac{|UV|}{2}$ instead of
$|SV|$:
$$d(f(U),f(V))\ge\sqrt{\frac{1}{4}|UV|\cdot
\xi^2\Big(\frac{|UV|}{2}\Big)+{\text{const}}}$$ and finally
$$\rho_f(t)\succeq \sqrt{t}\cdot
\xi(t)=\frac{t}{\sqrt{\ln t}\cdot\ln\ln t}.  $$\end{proof}

\section{Uniform embeddings of hyperbolic groups}\label{HypGroups}

Let us recall one of many equivalent definitions of a word
hyperbolic group (see~\cite{GroHypGroups, BH} for more details). A
geodesic triangle in a metric space is called {\it $\delta$-thin}
if any of its sides is contained in $\delta$-neighborhood of the
union of the other two sides. A geodesic metric space is called
{\it hyperbolic} if there exists $\delta>0$ such that any geodesic
triangle in this space is $\delta$-thin. A finitely generated
group is called {\it hyperbolic} (in a sense of M.~Gromov) if it
is hyperbolic as a metric space with word metric.

Recall that a hyperbolic group is {\it elementary} if it is
virtually cyclic. Any elementary hyperbolic group can be
quasi-isometrically embedded into a Hilbert space. Any
non-elementary hyperbolic group contains a quasi-isometric image
of the free group $F_2$~\cite{GroHypGroups}.

\begin{thm}
Any hyperbolic group admits a uniform embedding into Hilbert space
with the compression $$\rho(t)\succeq \frac{t}{\sqrt{\ln
t}\cdot\ln\ln t}.
$$

If $G$ is a non-elementary hyperbolic group, then for any uniform
embedding $f$ of $G$ in Hilbert space
$$\rho_f(t)\preceq\frac{t}{\sqrt{\ln t}}$$

\end{thm}

\begin{proof} We construct a uniform embedding as a composition of
uniform embeddings as follows. Any word hyperbolic group embeds
quasi-isometrically into a hyperbolic space $\mathbb H^n$ for some
$n$ \cite{bonk-shram}. The space $\mathbb H^n$ embeds
quasi-isometrically into a finite product $\prod_{i=1}^k T_i$ of
locally finite trees equipped with the
$l_1$-metric~\cite{dra-shro}. A finite product of metric spaces
with $l_1$-metric is quasi-isometric to the same product equipped
with $l_2$-metric. By Theorem~\ref{embed tree} each tree $T_i$
embeds into Hilbert space $\HH_i$ with compression function
$\succeq \frac{t}{\sqrt{\ln t}\cdot\ln\ln t}$. Therefore, the
product $\prod_{i=1}^k T_i$ embeds into the Hilbert space
$\HH=\prod_{i=1}^k \HH_i$ with compression function $\succeq
\frac{t}{\sqrt{\ln t}\cdot\ln\ln t}$.

Assuming that $G$ is a non-elementary hyperbolic group, one can
find an undistorted subgroup in $G$ isomorphic to the free group
on two generators (\cite{olsh}). Now the required upper bound on
the compression follows from corollary~\ref{compression tree}.
\end{proof}

\section{Uniform embeddings of CAT(0) cubical complexes}\label{CAT(0)}

In this section we use the method of~\cite{cam-nib} to generalize
the construction of uniform embedding from Theorem~\ref{embed
tree}. Recall that a {\it cubical complex} is a polyhedral complex
with geodesic metric in which each cell is isometric to a
Euclidean unit cube $[0,1]^n$ for some $n$ and the gluing maps are
isometries. We refer the reader to the book~\cite{BH} for
definition and properties of CAT(0) spaces.

Let us recall some notions that we will use in this Section
(see~\cite{cam-nib} and references therein for more details). A
midplane of a cube $[0,1]^n$ is the intersection of the cube with
a codimension $1$ hyperplane parallel to a coordinate hyperplane
and passing through the center of the cube. There are $n$
midplanes in every $n$-dimensional cube. Given an edge in a CAT(0)
cubical complex, there is a unique codimension $1$ {\it
hyperplane} in the complex which cuts the edge transversely in its
midpoint. This is obtained by developing the midplanes in the
cubes containing the edge. Any hyperplane in a CAT(0) cubical
complex separates it into two components.

Let $K$ be a CAT(0) cubical complex with metric $d$. We denote by
$K^{(1)}$ its $1$-skeleton (the union of all vertices and edges of
$K$) equipped with geodesic metric $d_1$. Then the $d_1$-distance
between two vertices of $K^{(1)}$ is equal to the number of
hyperplanes in $K$ separating them. In case if $K$ is a tree, any
hyperplane is a midpoint of some edge and the distance between two
vertices $U$ and $V$ is simply the number of midpoints separating
$U$ and $V$.

In case of a tree the hyperplanes which separate two vertices are
linearly ordered. In a higher dimensional cubical complex they are
not, but the notion of normal cube path introduced in~\cite{NR}
allows one to introduce a partial order and to modify the argument
from Theorem~\ref{embed tree}. A {\it cube path} is a sequence of
cubes $\CC = \{C_1, \dots ,C_n\}$, each of dimension at least 1,
such that the intersection $C_i\cap C_{i+1}$ is a single vertex
$V_i$ and $C_i$ is the (unique) cube of minimal dimension
containing $V_i$ and $V_{i+1}$. We denote by $V_0$ the vertex of
$C_1$ which is diagonally opposite to $V_1$, and by $V_n$ the
vertex of $C_n$ diagonally opposite to $V_{n-1}$. The vertices
$V_0$ and $V_n$ are called the initial vertex and the terminal
vertex respectively. A cube path is called {\it normal} if
$C_{i+1}\cap\st (C_i)=V_i$ for each $i$, where $\st (C_i)$ is the
union of all cubes which contain $C_i$ as a face (including $C_i$
itself).

For any two vertices $U$ and $V$ of a CAT(0) cubical complex there
is a unique normal cube path $\CC = \{C_1, \dots ,C_n\}$ from $U$
to $V$. A hyperplane separates $U$ and $V$ if and only if it
intersects (exactly) one of the cubes in the path $\CC$.

Suppose that $U$ and $V$ are adjacent vertices of a CAT(0) cubical
complex. Take any vertex $O$ and look at the normal cube paths
$\CC_U = \{u_1, \dots ,u_n\}$ from $U$ to $O$ and $\CC_V = \{v_1,
\dots ,v_n\}$ from $V$ to $O$. The following key property of
normal cube paths basically says that these two paths stay close
to each other: if a hyperplane $h$ separates both $U$ and $V$ from
$O$ and intersects a cube $u_i$, then it also intersects one of
the cubes $v_{i-1}$, $v_i$, $v_{i+1}$.

\begin{thm}\label{embed complex}
Let $K$ be a finite-dimensional, locally finite CAT(0) cubical
complex. There exists a uniform embedding $f$ of $K$ in Hilbert
space $\HH$ with the compression $$\rho_{f}(t) \succeq
\frac{t}{\sqrt{\ln t}\cdot\ln\ln t}$$
\end{thm}

\begin{proof} If a cubical complex $K$ is finite dimensional, then it is
quasi-isometric to its 1-skeleton $K^{(1)}$. We will construct a
uniform embedding $f$ of $K^{(1)}$ into Hilbert space with
$$\rho_{f}(t) \succeq \frac{t}{\sqrt{\ln t}\cdot\ln\ln t}$$

Our embedding is going to be linear on all edges of $K^{(1)}$,
therefore it is enough to define it on vertices. We fix a base
vertex $O\in K^{(1)}$. Our uniform embedding $f\colon K^{(1)}\to
\HH$ sends O to the origin of $\HH$. Let us fix a countable basis
in the space $\HH$. We would like to enumerate the elements of the
basis by all hyperplanes in $K$. If $h$ is a hyperplane in $K$, we
denote by $\vec{e}_h$ the corresponding unit vector of the basis.

Let $V$ be any vertex of $K^{(1)}$. We denote by $\|V\|$ the
number of cubes in the (unique) normal cube path $\CC_V = \{v_1,
\dots ,v_{\|V\|}\}$ from $V$ to $O$. We are going to define the
map $f$ on the vertex $V$ in such a way that the point $f(V)$ has
non-zero coordinates with respect to basic vectors corresponding
to the hyperplanes separating $V$ and $O$. To define the map $f$
we fix a function $\xi\colon \RR\to\RR$ called {\it weight
function} of the embedding $f$ and put $$
f(V)=\sum_{i=1}^{\|V\|}\sum_{h\cap v_i\ne\emptyset} \xi(i)\cdot
\vec{e}_{h}$$ for any vertex $V$ of $K^{(1)}$.

Let us define a function $N_V$ on the set of all hyperplanes:

$$
   N_V(h)=\begin{cases} i
                  &\text{if \;$h\cap v_i\ne\emptyset$}\\
                  0,
                  &\text{if \;$h\cap v_j=\emptyset$ for all $v_j\in\CC_V$}\end{cases}
$$
Assuming that $\xi(0)=0$ we can write $$ f(V)=\sum_{h}
\xi(N_V(h))\cdot \vec{e}_{h}$$

We are going to use essentially the same weight function we used
in Theorem~\ref{embed tree}. It is convenient here to assume that
$\xi(t)$ is non-decreasing. So, we use the formula
$$ \xi(t)=\begin{cases} \frac{\sqrt{t}}{\sqrt{\ln t}\cdot\ln\ln t}
                  &\text{if \;$t\ge M$}\\
                  0,
                  &\text{if \;$t<M$}\end{cases}
$$
where $M$ is some fixed positive integer such that the function
$\frac{\sqrt{t}}{\sqrt{\ln t}\cdot\ln\ln t}$ increases for $t>M$.

In order to show that $f$ is a uniform embedding we need to
estimate its dilatation $\delta_f$ from above and its compression
$\rho_f$ from below.

First we show that $\delta_f(t)\preceq  t$. To prove this we show
that there is a constant $C$ such that the expansion of any edge
of $K^{(1)}$ under the map $f$ does not exceed $C$. Let $U$ and
$V$ be adjacent vertices of $K^{(1)}$. We are going to find a
constant $C$ such that $d(f(U),f(V))\le C\cdot d(U,V)=C$. Clearly
$$[d(f(U),f(V))]^2=\sum_{h}\left(
\xi(N_U(h))-\xi(N_V(h))\right)^2$$ By the key property of normal
cube paths $N_U(h)$ is equal to one of $N_V(h)-1$, $N_V(h)$,
$N_V(h)+1$. Denote by $h_{UV}$ the hyperplane separating $U$ and
$V$. Without loss of generality we can assume that $N_U(h_{UV})=0$
and $N_V(h_{UV})=1$. Therefore
$$[d(f(U),f(V))]^2=\xi^2(1)+\sum_{i=1}^{\|U\|}\sum_{N_U(h)=i}\left(
\xi(i)-\xi(i\pm 1)\right)^2$$ Since the complex $K$ has finite
dimension $n$, the number of hyperplanes with $N_U(h)=i$ is at
most $n$. Thus
$$[d(f(U),f(V))]^2\le\sum_{i=1}^{\|U\|}2n\left(
\xi(i)-\xi(i+ 1)\right)^2\le 2n\sum_{i=1}^{\infty}\left(
\xi(i+1)-\xi(i)\right)^2$$ part (1) of Lemma~\ref{calc} implies
that the last sum is finite and therefore the distance
$d(f(U),f(V))$ is bounded by some constant $C$.

Now we estimate the compression function $\rho_f$. Let $U$ and $V$
be any vertices of the complex $K$. There are $d_1(U,V)$
hyperplanes in $K$ separating $U$ from $V$. Without loss of
generality we may assume that at least half of them separate $U$
from $O$; denote by $\hh$ the set of these hyperplanes. Clearly
$$[d(f(U),f(V))]^2\ge\sum_{h\in\hh}\xi^2(N_U(h))=\sum_i\sum_{N_U(h)=i, h\in\hh}\xi^2(i)$$
There are at least $\lfloor\frac{d_1(U,V)}{2}\rfloor$ hyperplanes
in the set $\hh$ (where $\lfloor t\rfloor$ denotes the largest
integer smaller than $t$). Since there are at most $n=\dim K$
hyperplanes with $N_U(h)=i$ for every $i$ and the function $\xi$
is non-decreasing, we have
$$\sum_i\sum_{N_U(h)=i, h\in\hh}\xi^2(i)\ge \sum_{i=1}^{\lfloor\frac{d_1(U,V)}{2n}\rfloor}n\cdot\xi^2(i)$$
Using part (2) of Lemma~\ref{calc} we estimate the last sum and
finally
$$[d(f(U),f(V))]^2\ge n\cdot
\frac{1}{2}\lfloor\frac{d_1(U,V)}{2n}\rfloor\cdot
\xi^2(\lfloor\frac{d_1(U,V)}{2n}\rfloor)+{\text{const}}$$
Therefore
$$\rho_f(t)\succeq \sqrt{t}\cdot
\xi(t)\sim\frac{t}{\sqrt{\ln t}\cdot\ln\ln t}.  $$
\end{proof}

\begin{thm}
If a group $G$ acts properly and cocompactly on a finite
dimensional CAT(0) cubical complex, then $G$ admits a uniform
embedding into Hilbert space with compression $$\rho(t) \succeq
\frac{t}{\sqrt{\ln t}\cdot\ln\ln t}.$$ Suppose that the orders of
finite subgroups of $G$ are uniformly bounded. Then, either $G$
can be quasi-isometrically embedded into Hilbert space or for any
uniform embedding $f$ into Hilbert space
$$\rho_f(t)\preceq\frac{t}{\sqrt{\ln t}}.$$
\end{thm}

\begin{proof}
Since the action is proper and cocompact, the group $G$ embeds
quasi-isometrically into a finite dimensional CAT(0) cubical
complex and the first statement of the theorem follows from
Theorem~\ref{embed complex}.

If $G$ acts properly and cocompactly on a finite dimensional
CAT(0) cubical complex and the orders of finite subgroups of $G$
are uniformly bounded, then either $G$ is virtually a finitely
generated abelian group or $G$ has a rank 2 free
subgroup~\cite{sag-wise}. In the former case $G$ embeds
quasi-isometrically into Hilbert space. In the latter case one can
choose a quasi-isometrically embedded rank $2$ free subgroup in
$G$ (the same proof works as in~\cite{sag-wise}). Application of
Corollary~\ref{compression tree} completes the proof.
\end{proof}

\end{document}